\input amstex
\documentstyle{amsppt}
\magnification=\magstep1

\pagewidth{6.5truein}
\pageheight{9.0truein}
\NoBlackBoxes

\input xy
\xyoption{matrix}\xyoption{arrow}\xyoption{curve}

\long\def\ignore#1{#1}

\def\edge{\ar@{-}}

\TagsOnRight

\def\la{\Lambda}
\def\laMod{\Lambda{\text{-}}\roman{Mod}}
\def\lamod{\Lambda{\text{-}}\roman{mod}}
\def\add{\operatorname{add}}
\def\lfindim{\operatorname{l\,fin\,dim}}
\def\lFindim{\operatorname{l\,Fin\,dim}}
\def\rfindim{\operatorname{r\,fin\,dim}}
\def\rFindim{\operatorname{r\,Fin\,dim}}
\def\Findim{\operatorname{Fin\,dim}}
\def\gldim{\operatorname{gl\,dim}}
\def\rep{\operatorname{rep}}
\def\lglrep{\operatorname{l\,gl\,rep}}
\def\rglrep{\operatorname{r\,gl\,rep}}

\def\lfinrep{\operatorname{l\,fin\,rep}}
\def\lFinrep{\operatorname{l\,Fin\,rep}}

\def\Hom{\operatorname{Hom}}
\def\Tor{\operatorname{Tor}}
\def\Ext{\operatorname{Ext}}
\def\pdim{\operatorname{p\,dim}}
\def\idim{\operatorname{i\,dim}}
\def\O{{\Cal O}}
\def\A{{\Cal A}}
\def\NN{{\Bbb N}}
\def\QQ{{\Bbb Q}}
\def\ZZ{{\Bbb Z}}
\def\opio{{\Cal O}/\pi{\Cal O}}
\def\bfb{{\bold B}}
\def\bfu{{\bold U}}
\def\bfv{{\bold V}}
\def\tr{\operatorname{tr}}

\def\rad{\operatorname{rad}}
\def\tlgam{\widetilde\Gamma}

\def\AnFu{{\bf 1}}
\def\BasI{{\bf 2}}
\def\Bas{{\bf 3}}
\def\CoFu{{\bf 4}}
\def\Fuj{{\bf 5}}
\def\FuSa{{\bf 6}}
\def\FuWa{{\bf 7}}
\def\GKK{{\bf 8}}
\def\domino{{\bf 9}}
\def\convenient{{\bf 10}}
\def\analyzing{{\bf 11}}
\def\IgZa{{\bf 12}}
\def\Jan{{\bf 13}}
\def\JaOd{{\bf 14}}
\def\Jat{{\bf 15}}
\def\Kap{{\bf 16}}
\def\KiKu{{\bf 17}}
\def\KKprivate{{\bf 18}}
\def\KKS{{\bf 19}}
\def\Rum{{\bf 20}}
\def\Tar{{\bf 21}}
\def\Tartwo{{\bf 22}}
\def\WiRo{{\bf 23}}
\def\Wil{{\bf 24}}
\def\Zak{{\bf 25}}

\topmatter
\title
 Repetitive resolutions over classical orders and finite dimensional
algebras
\endtitle
\rightheadtext{Repetitive resolutions}

\author K. R. Goodearl and B. Huisgen-Zimmermann\endauthor

\address Department of Mathematics, University of California, Santa Barbara,
CA 93106\endaddress
\email goodearl\@math.ucsb.edu\endemail

\address Department of Mathematics, University of California, Santa Barbara,
CA 93106\endaddress
\email birge\@math.ucsb.edu\endemail

\thanks The research of both authors was partially supported by National
Science Foundation grants.\endthanks

\abstract Repetitiveness in projective and injective resolutions and its influence on homological dimensions are studied. Some variations on the theme
of repetitiveness are introduced, and it is shown that the corresponding invariants lead to very
good -- and quite accessible -- upper bounds on various finitistic
dimensions in terms of individual modules. These invariants are the
`repetition index' and the `syzygy type' of a module $M$ over an artinian
ring $\Lambda$. The repetition index measures the degree of
repetitiveness among non-projective direct summands of the syzygies of $M$, while the syzygy type of $M$ measures the number of indecomposable modules
among direct summands of the syzygies of $M$. It is proved that if $T$ is a right
$\Lambda$-module which contains an isomorphic copy of $\Lambda/J(\Lambda)$, then the
left big finitistic dimension of $\Lambda$ is bounded above by the repetition index of $T$, which
in turn is bounded above by the syzygy type of $T$. 

The finite dimensional $K$-algebras
$\Lambda = {\Cal O}/\pi{\Cal O}$, where $\Cal O$ is a classical order over a discrete valuation ring
$D$ with uniformizing parameter $\pi$ and residue class field $K$, are investigated. It is proved that, if $\text{gl.dim.}\, {\Cal O} =d<\infty$, then the global repetition index of $\Lambda$ is $d-1$ and all finitely generated $\Lambda$-modules have finite
syzygy type. Various examples illustrating the results are presented.\endabstract

\endtopmatter

\document

\head I. Introduction\endhead

It has long been known that a great deal of homological information can be
gleaned from repetitions occurring in projective and injective resolutions. A
fundamental -- and ancient -- example is the $(\ZZ/p^n\ZZ)$-module
$\ZZ/p\ZZ$, whose successive syzygies alternate between $p\ZZ/p^n\ZZ$ and
$\ZZ/p\ZZ$; from this phenomenon one immediately deduces that the module has
infinite projective dimension. Numerous authors have systematized various
kinds of repetitive homological behavior, e.g., with the notion of an
`ultimately closed' projective resolution which is due to Jans \cite{\Jan},
or that of a resolution `with a strongly redundant image' which was
introduced by Colby and Fuller
\cite{\CoFu} (see also \cite{\FuSa} and \cite{\FuWa}). The
consequences of various types of repetitive homology have also been studied
by Igusa-Zacharia
\cite{\IgZa}, Kirkman-Kuzmanovich-Small \cite{\KKS}, and Wilson \cite{\Wil},
to name a few other instances among many.

In the first part of this note, we introduce some variations on the theme
of repetitiveness and show that the corresponding invariants lead to very
good -- and quite accessible -- upper bounds on various finitistic
dimensions in terms of individual modules. These invariants are the
`repetition index' and the `syzygy type' of a module $M$ over an artinian
ring $\la$. The repetition index, $\rep(M)$, measures the degree of
repetitiveness among non-projective direct summands of the syzygies of $M$;
to say that $\rep(M)=n<\infty$, means tightening somewhat the Colby-Fuller
concept of a `projective resolution with a strongly redundant image from
$n$'. The syzygy type of $M$ measures the number of indecomposable modules
among direct summands of the syzygies of $M$; in case the module $M$ is
finitely generated, finiteness of its syzygy type is equivalent to the
existence of an ultimately closed projective resolution in the sense of
Jans. We prove that if
$T$ is a right
$\la$-module which contains an isomorphic copy of $\la/\rad(\la)$, then the
left big finitistic dimension of $\la$ is bounded above by $\rep(T)$, which
in turn is bounded above by the syzygy type of $T$. On one hand, this point
of view leads to several methods for bounding finitistic dimensions; on the
other, it frequently allows us to reduce the computation of a host of
projective dimensions to the task of resolving a single module.

Our second goal here is to present a new class of algebras $\la$ for which
the global repetition index (the supremum of the repetition indices of all
$\la$-modules) is finite. These are the finite dimensional $K$-algebras
$\la=\opio$, where $\O$ is a classical order over a discrete valuation ring
$D$ with uniformizing parameter $\pi$ and residue class field $K$. The
homological properties of $\O$ are to a great extent determined by those of
$\la$, as was observed in \cite{\Tartwo}, \cite{\Jat}, and \cite{\GKK},
while the latter algebra is substantially easier to handle. We are indebted
to E\. Kirkman and J\. Kuzmanovich for bringing the crucial role played by
these algebras to our attention. From our
main result (Theorem 4.3) it follows that, if $\gldim\O =d<\infty$, then
$\lglrep\la =d-1$ and all finitely generated $\la$-modules have finite
syzygy type. To illustrate the transport of information between $\O$ and
$\la$, we present an example in which the left and right finitistic
dimensions of
$\O$ differ, and one in which $\O$ has infinite
global dimension while  $\la$ is Gorenstein; in particular, $\lfindim\O=
\rfindim\O$ in this latter example.

Finally, we wish to advertize the algebras $\opio$ within the
representation theory community, as constituting a very interesting class of
finite dimensional algebras deserving further study, both for their own sake
and for their impact on the theory of classical orders.

\subhead Terminology\endsubhead Given an artinian ring $\la$, we denote by
$\laMod$ and $\lamod$ the category of all left $\la$-modules and the full
subcategory of finitely generated left $\la$-modules, respectively.
Moreover, given a class $\A$ of left $\la$-modules, the notation $\add\A$
will stand for the full subcategory of $\laMod$ having as objects all
finite direct sums of copies of direct summands of objects from $\A$.

In addition to the global dimensions, we will consider the big and little
finitistic dimensions, $\lFindim\la$ and $\lfindim\la$; the former is the
supremum of all the \underbar{finite} projective dimensions attained on
objects in $\laMod$, while the latter is the analogous supremum in $\lamod$.

Throughout, $J$ will be the Jacobson radical of $\la$.

We next present a few known results which will be used in the sequel.

\proclaim{Proposition 1.1} If $\la$ is left noetherian, then $\lFindim\la
\leq \idim {}_\la\la$.\endproclaim

\demo{Proof} \cite{\Bas}, Proposition 4.3. \qed\enddemo

\proclaim{Proposition 1.2} If $\la$ is noetherian and $\idim {}_\la\la$ and
$\idim \la_\la$ are both finite, then $\lFindim\la= \rFindim\la=
\idim {}_\la\la= \idim \la_\la$.\endproclaim

\demo{Proof} \cite{\Zak}, Lemma A and \cite{\KKS}, Proposition 2.1.
\qed\enddemo

In case $\la$ is a finite dimensional algebra, there is an easier proof of
the above result, which, moreover, yields some additional information:

\proclaim{Proposition 1.3} Let $\la$ be a finite dimensional algebra over a
field $K$, with Jacobson radical $J$. If $\idim {}_\la\la$ and
$\idim \la_\la$ are both finite, then 
$$\align \lFindim\la &= \rFindim\la= \lfindim\la= \rfindim\la\\
 &= \idim {}_\la\la= \idim \la_\la = \pdim E(_\la(\la/J))=
\pdim E((\la/J)_\la). \endalign$$\endproclaim

\demo{Proof} If $\idim {}_\la\la= m$ and $\idim \la_\la=n$, then
$\lFindim\la
\le m$ and $\rFindim\la \le n$ by Proposition 1. The duality $\Hom_K(-,K)$
carries $_\la\la$ to the right minimal injective cogenerator
$E((\la/J)_\la)$, and so
$$m= \pdim E((\la/J)_\la) \le \rfindim\la \le \rFindim\la \le n$$
because $E((\la/J)_\la)$ is finitely generated. Similarly,
$$n= \pdim E(_\la(\la/J)) \le \lfindim\la \le \lFindim\la \le m.
\qed$$\enddemo

\subhead Graphing conventions\endsubhead For the convenience of the reader,
we include a sketch of our graphing conventions for modules over path
algebras modulo relations; these may differ from standard ones in that they
take the $J$-layering of the modules into account. Suppose that
$\la=K\Gamma/I$ is a path algebra modulo relations over a field $K$, where
$\Gamma$ is a quiver and $I$ an admissible ideal in $K\Gamma$. We will
label the vertices of $\Gamma$ by integers, and denote the primitive
idempotent in
$\la$ corresponding to a vertex $i$ by $e_i$; further, we will
often write $S_i$ for the simple left $\la$-module $\la e_i/Je_i$.

For example, suppose that $\Gamma$ is the quiver

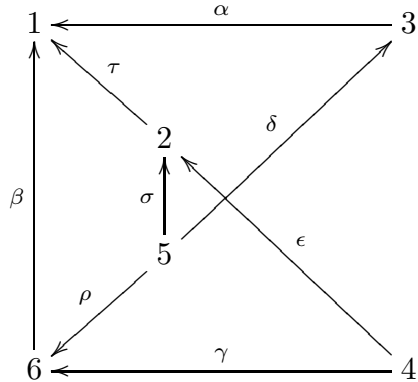
\ignore{
$$\xymatrixcolsep{2.5pc}
\xy\xymatrix{
1 &&&3 \ar[lll]_\alpha\\
 &2 \ar[ul]_\tau\\
 &5 \ar[u]^\sigma \ar[uurr]^\delta \ar[dl]_\rho\\
6 \ar[uuu]^\beta &&&4 \ar[lll]_\gamma \ar[uull]_\epsilon
}\endxy$$
}

\noindent That a left $\la$-module $M$ has graph

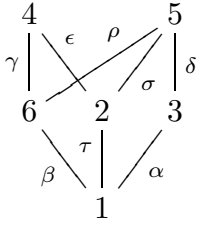
\ignore{
$$\xymatrixcolsep{1pc}\xymatrixrowsep{1.5pc}
\xy\xymatrix{
4 \edge[d]_\gamma \edge[dr]^(0.35)\epsilon &&5 \edge[dll]_(0.35)\rho
\edge[dl]^(0.6)\sigma \edge[d]^\delta\\ 
6 \edge[dr]_\beta &2 \edge[d]_(0.35)\tau &3 \edge[dl]^\alpha\\
 &1
}\endxy$$
}

\noindent is to communicate the following information:  First,
$$M/JM\cong S_4\oplus S_5; \quad JM/J^2M\cong S_6\oplus S_2\oplus S_3;
\quad J^2M\cong S_1; \quad J^3M=0.$$
Second, there exist elements $x=e_4x$ and $y=e_5y$ in $M\setminus JM$ such
that $\gamma x,\rho y\in JM\setminus J^2M$ and these last two elements
differ only by a nonzero scalar factor modulo $J^2M$; analogously,
$\epsilon x$ and $\sigma y$ belong to $JM\setminus J^2M$ and differ only
by a nonzero scalar modulo
$J^2M$. Finally, $\beta\gamma x$, $\tau\epsilon x$,
$\alpha\delta y$ are to differ only by nonzero scalar factors, that is,
$\la\beta\gamma x= \la\tau\epsilon x= \la\alpha\delta y= J^2M$ in this
case. (It follows from the information in the previous `layer' that
$\beta\rho y$ and $\tau\sigma y$ are automatically included in this list.)

Note that, in general, a module $M$ will have more than one graph, since
the latter may depend on the choice of the `reference elements' in
$M\setminus JM$ (the elements $x$ and $y$ in the example above). In case
there is just one arrow from a vertex $i$ to a vertex $j$, we will often
omit naming this arrow. For additional information, we refer to
\cite{\analyzing}.

Finally, observe that right $\la$-modules correspond to left modules over a
quotient of the path algebra of the `dual' of the quiver $\Gamma$, that is,
 of the quiver with the same vertices as $\Gamma$ but reversed arrows.
Thus in the graph of a \underbar{right} $\la$-module, an edge leading from a
vertex $j$ down to a vertex $i$ corresponds to an arrow
$i\rightarrow j$, rather than to an arrow $j\rightarrow i$.

\head 2. Syzygy type and repetition index\endhead

Throughout this section, let $\la$ be an artinian ring. There are analogs
of many of the results of this section for finitely generated
modules over a noetherian semiperfect ring -- see Remark 2.10.

\definition{Definitions 2.1} Let $M$ and $A$ be left $\la$-modules.

(1) Let $n$ be a nonnegative integer. We shall say that the minimal
projective resolution of $M$ is {\it repetitive at degree $n$} if there
exists a decomposition $\Omega^n(M)= P\oplus \bigoplus_{i\in I} A_i$ such
that $P$ is projective and each $A_i$ occurs as a direct summand of
infinitely many $\Omega^j(M)$. In case $\Omega^n(M)$ is finitely generated,
this is the same as to say that each non-projective indecomposable direct
summand of $\Omega^n(M)$ occurs as a direct summand of infinitely many
$\Omega^j(M)$.

(2) The {\it repetition index} of $M$, denoted $\rep(M)$, is the least
nonnegative integer $k$ such that the minimal projective resolution of $M$
is repetitive at degree $k$ (if such a
$k$ exists) or
$\infty$ (otherwise). The corresponding globalized indices are
$$\align \lglrep(\la) &= \sup\{ \rep(M) \mid M\in\laMod \}\\
\lFinrep(\la) &= \sup\{ \rep(M) \mid M\in\laMod \text{\ and\ } \rep(M)
<\infty\}\\
\lfinrep(\la) &= \sup\{ \rep(M) \mid M\in\lamod \text{\ and\ } \rep(M)
<\infty\}, \endalign$$
which we shall call the {\it left global repetition index}, the {\it left big
finitistic repetition index}, and the {\it left little finitistic repetition
index} of $\la$, respectively.

(3) Next, we set $\sigma_M(A) =-1$ if $A$ is not isomorphic to a direct
summand of any syzygy $\Omega^j(M)$, $j\ge 0$, and 
$$\sigma_M(A)= \sup\{ j\ge 0\mid A\cong \text{\ a direct summand of\ }
\Omega^j(M) \}$$
otherwise. We call $\sigma_M(A)$ the {\it (homological) contingency of $A$ to
$M$}.

Clearly $\rep(M)$ equals the minimum of those $k\ge 0$ such that
$\Omega^k(M)$ is a direct sum of a projective module and modules with
infinite contingency to $M$.

(4) The {\it syzygy category of $M$} is the full subcategory $\add\bigl( \{
\Omega^j(M) \mid j\ge 0\} \bigr)$ of $\laMod$. If this category is a
Krull-Schmidt category of finite representation type, then we say that $M$
has {\it finite syzygy type\/}, and we define the {\it syzygy type} of
$M$ to be the number of isomorphism types of nonzero indecomposable objects in
the syzygy category. Otherwise, we say that the syzygy type of $M$ is
$\infty$.

(5) For any nonnegative integer $\tau$, let $\Omega^\tau(\laMod)$ denote the
full subcategory 
$$\add\bigl( \{ \Omega^j(M) \mid j\ge \tau;\ M\in \laMod \} \bigr)$$
of $\laMod$. As in (4), if $\Omega^\tau(\laMod)$ is a Krull-Schmidt
category of finite representation type, then we say that $\laMod$ has {\it
finite syzygy type from degree $\tau$}, and we define the {\it syzygy type
of $\laMod$ from degree $\tau$} to be the number of isomorphism types of
nonzero indecomposable objects in the category $\Omega^\tau(\laMod)$.
\enddefinition

\definition{Comments 2.2} (a) The earliest consideration of these
conditions goes back to Jans \cite{\Jan}, who said that a 
$\la$-module $M$ has an {\it ultimately closed projective resolution} in
case some syzygy $\Omega^n(M)$, $n\ge1$, has a decomposition $\Omega^n(M)
=\bigoplus_{i\in I} A_i$ such that each $A_i$ already occurred as a
direct summand of a previous syzygy $\Omega^{j(i)}(M)$ with $j(i) <n$. It
is easily seen that, for a finitely generated module $M$, this is the same
as to say that $M$ has finite syzygy type. Subsequently, Colby-Fuller
\cite{\CoFu} and Fuller-Wang
\cite{\FuWa} relaxed this condition slightly as follows: $M$ is said to have
a projective resolution with a {\it strongly redundant image from an integer
$n\ge 1$} in case the
$n$-th syzygy has a decomposition $\Omega^n(M) = \bigoplus_{i\in I} A_i$
such that each $A_i$ is a direct summand of a \underbar{later} syzygy
$\Omega^{j(i)}(M)$ with $j(i) >n$. In case $\Omega^n(M)$ has a
 decomposition into modules with local endomorphism rings, this latter
condition clearly implies $\rep(M)\le n$. Note that strict inequality may
occur, because we exclude projective direct summands of the $\Omega^j(M)$'s
from consideration in $\rep(M)$. We point to \cite{\CoFu} for an example of a
projective resolution with a strongly redundant image
which fails to be ultimately closed.

(b) Observe that the restriction to non-projective direct summands in our
definition of `repetition index' allows us to compute the repetition index
of a finitely generated module $M$ by looking at \underbar{any} projective
resolution, not just a minimal one.

(c) If $\laMod$ has finite syzygy type from some degree, then each left
$\la$-module has finite repetition index. In fact, $\lglrep(\la) <\infty$ in
that case -- see Lemma 2.4(b).

(d) Note that $\pdim(M)= \rep(M)$ for any module $M$ with finite projective
dimension. Hence, we have the inequalities
$\lfindim\la
\le
\lfinrep\la
\le
\lglrep\la$ and $\lFindim\la \le \lFinrep\la \le \lglrep\la$.
Strict inequalities may occur in all places; for an instance where even
$\lFindim\la < \lfinrep\la$, see Example 2.3(d).
\enddefinition

\definition{Examples 2.3} (a) The module category $\laMod$ of any finite
dimensional monomial relation algebra
$\la= K\Gamma/I$ has finite syzygy type from degree 2, by \cite{\domino},
Theorem A; in fact, the syzygy type of $\laMod$ from degree 2 is bounded
above by the number of paths in $K\Gamma$ whose residue classes are
nontrivial in
$\la$.

(b) If $\gldim\la =m<\infty$, then $\la$ has finite syzygy type from degree
$m$.

(c) If $\la= K[X,Y]/(X^2,Y^2)$ for a field $K$, then $\la$ has infinite
repetition index and, a fortiori, infinite syzygy type; in fact, $\rep(\la/J)
=\infty$. On the other hand, $\Findim\la =\rep(E(\la/J)) =0$.

(d) Let $\la= K\Gamma/I$, where $\Gamma$ is the quiver

\ignore{
$$\xymatrixcolsep{3pc}
\xy\xymatrix{
1 \ar[r]^\alpha &2 \ar[r]^\beta &3 \ar@(ur,dr)^\gamma
}\endxy$$
}

\noindent and $I$ is generated by all paths of length 3. It is easy
to see that $\lfindim\la =\lFindim\la =0$ in this case (see \cite{\BasI},
Theorem 6.3 and \cite{\convenient}, Corollary 8). On the other hand,
$\lfinrep(\la) \ge \rep(S_1)= 2$, since $\Omega^1(S_1)$ has graph 
\ignore{
$\vcenter{
\xymatrixrowsep{1pc}
\xy\xymatrix{
2 \edge[d]\\ 3
}\endxy}$,
}
while $\Omega^2(S_1)
\cong \Omega^{2n}(S_1) \cong S_3$ for $n\ge 1$, and $\Omega^3(S_1)\cong
\Omega^{2n+1}(S_1)$ has graph 
\ignore{
$\vcenter{
\xymatrixrowsep{1pc}
\xy\xymatrix{
3 \edge[d]\\ 3
}\endxy}$
}
for $n\ge 1$. Thus indeed, $\lFindim\la < \lglrep\la$.

Finally, observe that the syzygy category of $S_1$ is the additively closed
category generated by the indecomposable left $\la$-modules $\Omega^i(S_1)$,
$i=0,1,2,3$, in this example. \qed
\enddefinition

The principal connections among syzygy type, repetition index, and
finitistic dimensions are as follows. Note that, in essence, most of the
arguments are familiar -- see, e.g., \cite{\IgZa}, \cite{\FuWa}. In fact, for
$T=\la/J$, where $\la$ is an artin algebra, the first statement of
Theorem 2.6 is implicit in \cite{\IgZa}.

\proclaim{Lemma 2.4} {\rm (a)} The repetition index of any $\la$-module $M$
is bounded above by the syzygy type of $M$.

{\rm (b)} If $\laMod$ has syzygy type $s$ from degree $\tau$,
then $\lglrep(\la) \le \tau+s$.\endproclaim

\demo{Proof} (a) Assume that the syzygy type of $M$ is $s<\infty$. If $s=0$,
then $M=0$ and hence $\rep(M)=0$. 

Now assume that $s>0$, and let $B=B_s$ be an indecomposable direct summand
of $\Omega^s(M)$. Since $M$ has finite syzygy type, $\Omega^{s-1}(M)$ is a
direct sum of indecomposable modules $I_k$, whence $\Omega^s(M) \cong
\bigoplus_k \Omega^1(I_k)$, and so $B_s$ is isomorphic to a direct summand
of some $\Omega^1(I_k)$. (Here we are using the assumption that the syzygy
category of $M$ is a Krull-Schmidt category.) Label this
$I_k$ as
$B_{s-1}$. Continuing by induction, we obtain indecomposable modules $B_s,
B_{s-1},\dots, B_0$ such that each $B_j$ is a direct summand of
$\Omega^j(M)$ and $B_j$ is isomorphic to a direct summand of
$\Omega^1(B_{j-1})$ for $j>0$.

By definition of $s$, there exist indices $p,q$ with $0\le p<q\le s$ such
that $B_p\cong B_q$. Thus $B_p$ is isomorphic to a direct summand of
$\Omega^{q-p}(B_p)$, and hence isomorphic to a direct summand of
$\Omega^{m(q-p)}(B_p)$ for all $m\ge 0$. Since $B$ is isomorphic to a
direct summand of $\Omega^{s-p}(B_p)$ and $B_p$ is isomorphic to a direct
summand of $\Omega^p(M)$, it follows that $B$ is isomorphic to a direct
summand of $\Omega^{s+m(q-p)}(M)$ for all $m\ge 0$. Therefore $\rep(M)\le
s$, because $B$ was an arbitrary indecomposable direct summand of
$\Omega^s(M)$.

(b) It follows from our assumptions that for any $M\in \laMod$, the module
$\Omega^\tau(M)$ has syzygy type at most $s$. Thus $\rep(\Omega^\tau(M))
\le s$ by part (a), and therefore $\rep(M)\le \tau+s$.
\qed\enddemo 

\proclaim{Observation 2.5} Whenever $T$ is a right $\la$-module such that
$\la/J$, viewed as a right $\la$-module, embeds into $T$, 
$$\lFindim\la \le \rep(T).$$\endproclaim

\demo{Proof} Assume that $\rep(T) =r <\infty$, and write $\Omega^r(T)=
P\oplus \bigoplus_{i\in I} A_i$ where $P$ is projective and $\sigma_T(A_i)
=\infty$ for all $i$. Consider a left
$\la$-module
$M$ with
$\pdim M=d<\infty$, and note that
$\Tor_d^\la(\la/J,M)\ne 0$. By assumption, there exists a short exact
sequence of right
$\la$-modules
$0\rightarrow
\la/J\rightarrow T\rightarrow C\rightarrow 0$. Since
$\Tor^\la_{d+1}(C,M)=0$, we must have
$\Tor^\la_d(T,M)\ne 0$.

If $d>r$, then $\Tor^\la_{d-r}(\Omega^r(T),M)\ne 0$, and so
$\Tor^\la_{d-r}(A_i,M)\ne 0$ for some $i$. Since $\sigma_T(A_i) =\infty$, it
follows that
$\Tor^\la_{d-r}(\Omega^j(T),M)\ne 0$ for infinitely many $j$, and
consequently
$\Tor^\la_{d-r+j}(T,M)\ne 0$ for infinitely many $j$. However, this is
impossible for $j>r$. Therefore $d\le r$.

Since $M$ was an arbitrary left $\la$-module of finite projective dimension,
this proves that $\lFindim\la \le r$. \qed\enddemo

\proclaim{Theorem 2.6} Suppose there exists a right $\la$-module $T$ of
finite syzygy type $s$ such that $\la/J$, viewed as a right $\la$-module,
embeds into $T$. Moreover, let $A_1,\dots,A_s$ be representatives for the
isomorphism types of the indecomposable objects in the syzygy category of $T$.
Then

{\rm (a)} $\lFindim\la \le \rep(T) \le s$.

{\rm (b)} If $M$ is a left $\la$-module of finite projective dimension and
$$\mu= \max\{ \sigma_T(A_i) \mid 1\le i\le s \text{\ and\ } \Tor_1^\la(
A_i,M)\ne 0 \},$$
then $\pdim(M) = \mu+1$. (Let $\mu= -1$ if the above set is
empty.)\endproclaim

\demo{Proof} (a) follows from Lemma 2.4 and Observation 2.5.

(b) Set $d=\pdim M<\infty$, and note that
$\Tor^\la_d(T,M)\ne 0$, as in the previous proof.

Consider one of the $A_i$ for which $\Tor^\la_1(A_i,M)\ne 0$, and choose an
index
$j$ such that $A_i$ is isomorphic to a direct summand of $\Omega^j(T)$. Then
$\Tor^\la_{j+1}(T,M)\cong \Tor^\la_1(\Omega^j(T),M)\ne 0$, and so $j+1\le
d$. Thus $\sigma_T(A_i)+1\le d$ for all $A_i$ such that
$\Tor^\la_1(A_i,M)\ne 0$, that is, $\mu+1\le d$.

It remains to prove the reverse inequality. If $d=0$, then
$\Tor^\la_1(-,M)\equiv 0$, whence $\mu=-1$ and $d=\mu+1$. If $d>0$, then
$\Tor^\la_1(\Omega^{d-1}(T),M)\cong \Tor^\la_d(T,M)\ne 0$, and so
$\Tor^\la_1(A_i,M)\ne 0$ for some $A_i$ which is isomorphic to a direct
summand of $\Omega^{d-1}(T)$. In this case, $d-1\le \sigma_T(A_i)\le \mu$,
and therefore $d\le \mu+1$.
\qed\enddemo

\proclaim{Corollary 2.7} Let $T$ and $A_1,\dots,A_s$ be right
$\la$-modules as in the hypotheses of Theorem 2.6. Then
$$\lFindim\la \le 1+ \max_{1\le i\le s} \{ \sigma_T(A_i) \mid \ A_i \
\text{non-projective and} \ \sigma_T(A_i) <\infty \}. \qed$$
\endproclaim

Of course, the module $T$ in (2.5--2.7) may be chosen to be either the right
module $\la/J$ or the minimal injective cogenerator $E$ for the category
of right $\la$-modules. In \cite{\FuWa}, certain related bounds on the
finitistic dimensions of noetherian rings are obtained. For the case where
the base ring is an artin  algebra, our Observation 2.5 provides a mild
tightening of these bounds: For Theorem 3 of \cite{\FuWa}, take $T=E$; for
Theorem 9, take $T= \la/J$. 

For reasons that are not transparent as yet, the choice $T=E$ appears to
systematically yield better bounds on $\lFindim\la$. In particular, while it
is easy to construct algebras with $\rep (\la/J) =\infty$ (see Example
2.3(c) above), examples where $\rep(E) =\infty$ are not immediate.

The following example will illustrate the methods just introduced; further,
it will be relevant to Example 4.6.

\definition{Example 2.8} Let $\la= K\Gamma/I$ be a binomial relation algebra
(i.e., $I$ can be generated by paths in $K\Gamma$ and by differences $p-kq$
with $k\in K$, where $p$ and $q$ are paths) with six vertices such that the
indecomposable projective right $\la$-modules $P_1,\dots,P_6$ have the
following graphs:

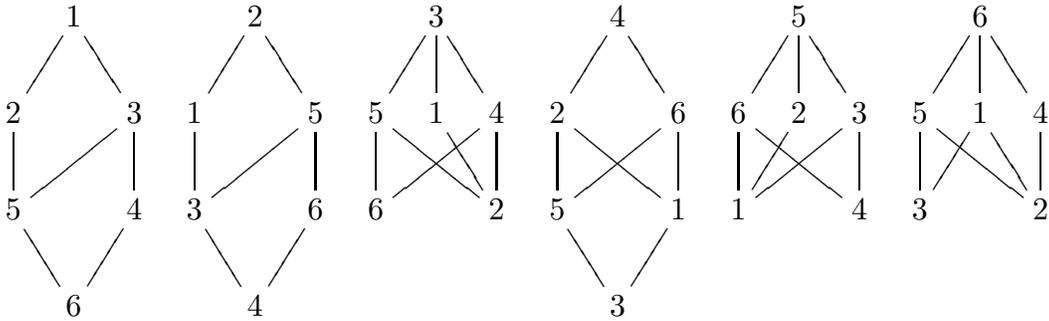
\ignore{
$$\xymatrixcolsep{0.67pc}\xymatrixrowsep{1.5pc}
\xy\xymatrix{
 &1 \edge[dl] \edge[dr] &&&2 \edge[dl] \edge[dr] &&&3 \edge[dl] \edge[d]
\edge[dr] &&&4 \edge[dl] \edge[dr] &&&5 \edge[dl] \edge[d] \edge[dr] &&&6
\edge[dl] \edge[d] \edge[dr]\\
2 \edge[d] &&3 \edge[dll] \edge[d] &1 \edge[d] &&5 \edge[dll] \edge[d] &5
\edge[d] \edge[drr] &1 \edge[dr] &4 \edge[dll] \edge[d] &2 \edge[d]
\edge[drr] &&6 \edge[dll] \edge[d] &6 \edge[d] \edge[drr] &2 \edge[dl] &3
\edge[dll] \edge[d] &5 \edge[d] \edge[drr] &1 \edge[dl] \edge[dr] &4
\edge[d]\\
5 \edge[dr] &&4 \edge[dl] &3 \edge[dr] &&6 \edge[dl] &6 &&2 &5 \edge[dr]
&&1 \edge[dl] &1 &&4 &3 &&2\\
 &6 &&&4 &&& &&&3
}\endxy$$
}

\noindent Then the indecomposable injective right $\la$-modules $E_i$ are
$E_3= P_4$, $E_4= P_2$, $E_6= P_1$, while $E_1$, $E_2$, $E_5$ have
graphs

\ignore{
$$\xymatrixcolsep{1pc}\xymatrixrowsep{1.5pc}
\xy\xymatrix{
4 \edge[d] \edge[dr] &&5 \edge[dll] \edge[dl] \edge[d] &&3 \edge[d]
\edge[dr] \edge[drr] &&6 \edge[dll] \edge[dl] \edge[d] &&1 \edge[d]
\edge[dr] &&4 \edge[dl] \edge[d]\\
6 \edge[dr] &2 \edge[d] &3 \edge[dl] &&5 \edge[dr] &1 \edge[d] &4
\edge[dl] &&3 \edge[dr] &2 \edge[d] &6 \edge[dl]\\
 &1 &&&&2 &&&&5
}\endxy$$
}

\noindent respectively.  Set $E= E_1 \oplus \cdots \oplus E_6$. Clearly
$\rep E_3= \rep E_4= \rep E_6 =0$. To compute $\rep E$, we determine the
minimal projective resolutions of $E_1$, $E_2$, $E_5$:

\ignore{
$$\xymatrixcolsep{1pc}\xymatrixrowsep{1.5pc}
\xy\xymatrix{
 &\Omega^1(E_1) &&&&\Omega^2(E_1) &&&&\Omega^3(E_1)\\
 &6 \edge[dl] \edge[d] \edge[dr] &2 \edge[dl] \edge[d] &&&5 \edge[dl]
\edge[d] \edge[dr] &1 \edge[dl] \edge[d] &&&3 \edge[dl] \edge[d]
\edge[dr] &2 \edge[dll] \edge[d]\\
4 &5 \edge[d] &1 \edge[dl] &&6 \edge[dr] &3 \edge[d] &2 &&1 &4 \edge[d]
&5 \edge[dl]\\
 &3 &&&&4 &&&&6
}\endxy$$
}

\noindent and $\Omega^4(E_1) \cong \Omega^2(E_1)$. Moreover, $\Omega^1(E_2)$
has graph

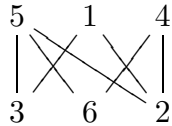
\ignore{
$$\xymatrixcolsep{1pc}\xymatrixrowsep{1.5pc}
\xy\xymatrix{
5 \edge[d] \edge[dr] \edge[drr] &1 \edge[dl] \edge[dr] &4 \edge[dl]
\edge[d]\\
3 &6 &2
}\endxy$$
}

\noindent and $\Omega^2(E_2) \cong \Omega^1(E_1) \oplus \Omega^3(E_1)$.
Finally, $\Omega^1(E_5)$ has graph

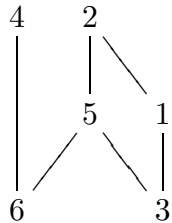
\ignore{
$$\xymatrixcolsep{1pc}\xymatrixrowsep{1.5pc}
\xy\xymatrix{
4 \edge[dd] &2 \edge[d] \edge[dr]\\
 &5 \edge[dl] \edge[dr] &1 \edge[d]\\
6 &&3
}\endxy$$
}

\noindent and $\Omega^2(E_5) \cong \Omega^1(E_1)$. It follows that $\rep E=
\rep(E_1\oplus E_2\oplus E_5) =3$, and hence Observation 2.5 yields
$\lFindim\la
\le 3$. We actually obtain equality, $\lFindim\la= \lfindim\la =3$, since the
indecomposable injective left $\la$-module $E(\la e_3/Je_3)$ has projective
dimension 3. Note that, on the other hand, $\rFindim\la =0$, since the left
socle of
$\la$ contains a copy of $\la/J$ (see \cite{\BasI}, Theorem 6.3 and
\cite{\convenient}, Corollary 8).
\qed
\enddefinition

Often the following variant of Theorem 2.6 provides better bounds:

\proclaim{Corollary 2.9} Suppose there exists a right $\la$-module $T$ of
finite syzygy type such that $\la/J$, viewed as a right $\la$-module,
embeds into $T$. For $m\ge 0$, let $s_m$ be the syzygy type of $\Omega^m(T)$.
Then $\lFindim\la \le s_m+m$ for all $m$. \endproclaim

\demo{Proof} By Lemma 2.4, $\rep(\Omega^m(T))\le s_m$, and so $\rep(T)\le
s_m+m$. The corollary thus follows immediately from part (a) of Theorem 2.6.
\qed\enddemo

\definition{Remark 2.10} Since the definitions and proofs in this section
rely mainly on the existence and behavior of syzygies, the results carry
over to the case when $\la$ is only assumed to be a noetherian semiperfect
ring, provided we restrict attention to finitely generated modules. (Recall
that semiperfect rings are precisely the rings over which every finitely
generated module has a projective cover.) In particular, the argument of
Lemma 2.4(a) shows that the repetition index of any finitely generated
$\la$-module is bounded above by its syzygy type, and the argument of
Observation 2.5 shows that if $T$ is a finitely generated right
$\la$-module containing an isomorphic copy of $(\la/J)_{\la/J}$, then
$$\lfindim\la \le \rep(T).$$
Similarly, there are analogs of Theorem 2.6(b) (giving a formula for the
projective dimensions of finitely generated $\la$-modules with finite
projective dimension) and of Corollaries 2.7, 2.9 (giving additional upper
bounds for $\lfindim\la$). We leave the details to the interested reader.
\enddefinition

\head 3. Computation of individual projective dimensions in the case
of finite syzygy type\endhead

We now specialize to the case where $\la$ is a finite dimensional algebra
over a field
$K$. Suppose there exists a finitely generated right $\la$-module $T$ of
finite syzygy type $s$ such that $(\la/J)_\la$ embeds into $T$, and let
$A_1,\dots,A_s$ be representatives for the isomorphism types of the
indecomposable objects in the syzygy category of $T$.

Guiding idea: If one knows the $A_i$ and their first syzygies, the
computation of projective dimensions of left $\la$-modules $M$ can be reduced
to the calculation of the
$K$-dimensions of the tensor products $A_i\otimes_\la M$.

Indeed, by the choice of the $A_i$, there exists an $s\times s$ matrix $\bfb=
(b_{ij})$ of nonnegative integers such that
$$\Omega^1(A_i)\cong \bigoplus_{j=1}^s A_j^{b_{ij}}$$
for all $i$. Moreover, let $e_1,\dots, e_n$ be a full set of orthogonal
primitive idempotents in $\la$, and let $P_i$ denote the projective cover of
$A_i$. Write
$$P_i\cong \bigoplus_{l=1}^n (e_l\la)^{p_{il}}$$
with $p_{il}\ge0$ for $1\le i\le s$.

Observe that since the entries of $\bfb$ are nonnegative integers, the
product of $\bfb$ with any column vector of cardinal numbers is defined. (We
use the convention that $0\cdot\alpha =0$.)

\proclaim{Proposition 3.1} Let $M$ be a left
$\la$-module, set $\tau_i= \dim_K \Tor^\la_1(A_i,M)$ for $i=1,\dots,s$, and
set $\tau =\tau(M)= (\tau_1,\tau_2,\dots,\tau_s)^{\tr}$.

{\rm (a)} If $M$ has finite projective dimension, then $\bfb^m\tau=0$
for some $m\ge 0$, and the least such $m$ equals $\pdim M$.

{\rm (b)} Suppose that $T= (\la/J)_\la$. Then $\pdim M<\infty$ if and only
if 
$\bfb^m\tau=0$ for some $m$.

{\rm (c)} If $M$ is finitely generated, the integers $\tau_i$, $1\le i\le
s$, can be computed as follows:
$$\tau_i= \dim_K
(A_i\otimes_\la M) +\sum_{j=1}^s b_{ij}
\dim_K(A_j\otimes_\la M) -\sum_{l=1}^n p_{il} \dim_K (e_lM).$$
\endproclaim

\demo{Proof} For $i=1,\dots,s$, we have
$$\Tor^\la_1(A_i,\Omega^1(M)) \cong \Tor^\la_2(A_i,M) \cong
\Tor^\la_1(\Omega^1(A_i),M) \cong \bigoplus_{j=1}^s
\Tor^\la_1(A_j,M)^{b_{ij}},$$
whence $\dim_K \Tor^\la_1(A_i, \Omega^1(M))= \sum_{j=1}^s b_{ij} \tau_j$.
This shows that $\tau(\Omega^1(M))= \bfb \tau(M)$. It follows by induction
that
$$\bfb^m\tau(M)= \tau(\Omega^m(M))$$
for all $m\ge 0$.

(a) Suppose that $\pdim M=m<\infty$. Then $\Omega^m(M)$ is projective, and
consequently $\Tor^\la_1(A_i,\Omega^m(M))=0$ for all $i$.
Hence, $\bfb^m\tau = \tau(\Omega^m(M)) =0$. We must also show that if $m>0$,
then $\bfb^{m-1}\tau \ne 0$. As in the proof of Observation 2.5,
$\Tor^\la_m(T,M) \ne 0$, whence $\Tor^\la_1(T,\Omega^{m-1}(M)) \ne 0$. Since
$T$ is a direct sum of copies of the $A_i$, it follows that $\Tor^\la_1(A_i,
\Omega^{m-1}(M)) \ne 0$ for some $i$, and therefore $\bfb^{m-1}\tau=
\tau(\Omega^{m-1}(M)) \ne 0$ as desired.

(b) Suppose that $\bfb^m\tau =0$ for some $m$. Then
$\tau(\Omega^m(M))=0$, and so we have $\Tor^\la_1(A_i, \Omega^m(M))=0$ for all
$i$. Consequently, $\Tor^\la_1(\la/J, \Omega^m(M))=0$, which implies that
$\Omega^m(M)$ is projective, and thus $\pdim M\le m$.

(c) From the resolutions $0\rightarrow \Omega^1(A_i)\rightarrow
P_i\rightarrow A_i\rightarrow 0$, we obtain long exact sequences
$$0\rightarrow \Tor^\la_1(A_i,M) \rightarrow \Omega^1(A_i)\otimes_\la M
\rightarrow P_i\otimes_\la M \rightarrow A_i\otimes_\la M \rightarrow 0,$$
and consequently
$$\tau_i= \dim_K(A_i\otimes_\la M) +\dim_K(\Omega^1(A_i)\otimes_\la M)
-\dim_K(P_i\otimes_\la M).$$
Since $\Omega^1(A_i)\otimes_\la M\cong \bigoplus_{j=1}^s (A_j\otimes_\la
M)^{b_{ij}}$ and $P_i\otimes_\la M\cong \bigoplus_{l=1}^n (e_lM)^{p_{il}}$,
the desired formula is clear. \qed\enddemo

\proclaim{Corollary 3.2} If $d$ is the least nonnegative integer with the
property that
$\bfb^d\in M_s(\QQ)\bfb^{d+1}$, then $\lFindim\la \le d$.\endproclaim

\demo{Proof} Observe that $\bfb^d\in M_s(\QQ)\bfb^{d+k}$ for all $k>0$.
Hence, there exist positive integers $z_k$ and nonnegative integer matrices
$\bfu_k,\bfv_k$ such that $z_k\bfb^d= (\bfu_k-\bfv_k)\bfb^{d+k}$.

Given any left $\la$-module $M$ with $\pdim M= m<\infty$, we have
$\bfb^m\tau(M) =0$ by Proposition 3.1. If $m>d$, multiply both sides of the
equation $z_{m-d}\bfb^d +\bfv_{m-d}\bfb^m =\bfu_{m-d}\bfb^m$ on the right by
$\tau(M)$ to obtain $z_{m-d}\bfb^d\tau(M) =0$. But then $\bfb^d\tau(M) =0$,
contradicting Proposition 3.1. Therefore $\pdim M\le d$. \qed\enddemo

We conclude this section with an example illustrating the use of Proposition
3.1.

\definition{Example 3.3} Let $\la= K\Gamma/I$ be the monomial relation
algebra with quiver

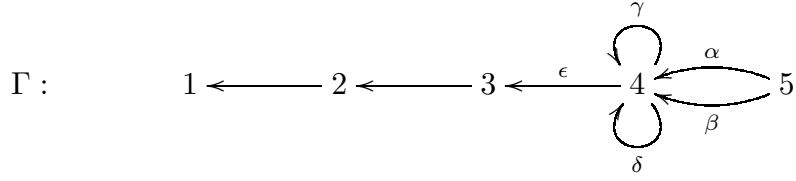
\ignore{
$$\xymatrixcolsep{3pc}
\xy\xymatrix{
\Gamma: &1 &2 \ar[l] &3 \ar[l] &4 \ar[l]_\epsilon \ar@(ur,ul)_\gamma
\ar@(dr,dl)^\delta &5 \ar@/_/[l]_\alpha \ar@/^/[l]^\beta
}\endxy$$
}

\noindent such that the graphs of the indecomposable projective left
$\la$-modules are

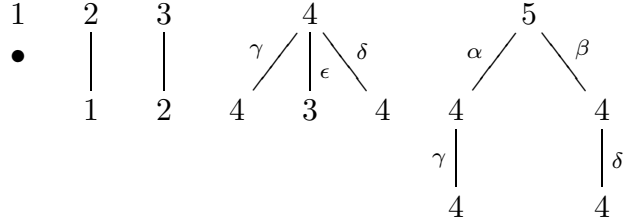
\ignore{
$$\xymatrixcolsep{1pc}\xymatrixrowsep{1.5pc}
\xy\xymatrix{
1 \save+<0ex,-3ex>\drop{\bullet}\restore &2 \edge[d] &3 \edge[d] &&4
\edge[dl]_\gamma \edge[d]^(0.6)\epsilon \edge[dr]^\delta &&&5
\edge[dl]_\alpha
\edge[dr]^\beta\\
 &1 &2 &4 &3 &4 &4 \edge[d]_\gamma &&4 \edge[d]^\delta\\
 &&&&&&4 &&4
}\endxy$$
}

\noindent Then the indecomposable projective right $\la$-modules have graphs

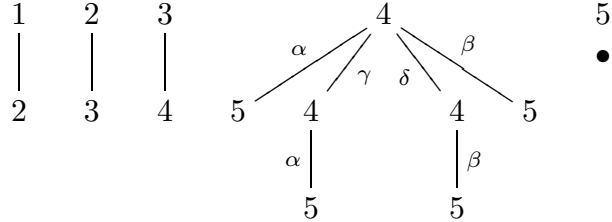
\ignore{
$$\xymatrixcolsep{1pc}\xymatrixrowsep{1.5pc}
\xy\xymatrix{
1 \edge[d] &2 \edge[d] &3 \edge[d] &&&4 \edge[dll]_\alpha
\edge[dl]^\gamma \edge[dr]_\delta \edge[drr]^\beta &&&5
\save+<0ex,-3ex>\drop{\bullet}\restore\\
2 &3 &4 &5 &4 \edge[d]_\alpha &&4 \edge[d]^\beta &5\\
 &&&&5 &&5
}\endxy$$
}

\noindent Let $e_1,\dots,e_5$ be the primitive idempotents of $\la$
corresponding to the vertices of $\Gamma$, set $S_i= e_i\la/e_iJ$, and let
$T$ be the right $\la$-module $\la/J\cong \bigoplus_{i=1}^5 S_i$. One
computes that
$$\Omega^1(S_1)\cong S_2; \qquad \Omega^1(S_2)\cong S_3; \qquad
\Omega^1(S_3)\cong S_4; \qquad \Omega^1(S_5) =0;$$
while
$$\Omega^1(S_4)\cong S_5^2 \oplus
\ignore{\vcenter{
\xymatrixrowsep{1pc}
\xy\xymatrix{
4\edge[d]_\alpha\\
5}\endxy}} 
\oplus
\ignore{\vcenter{
\xymatrixrowsep{1pc}
\xy\xymatrix{
4\edge[d]_\beta\\
5}\endxy}}
\qquad \text{and} \qquad
\Omega^2(S_4)\cong S_5^2 \oplus
\left(
\ignore{\vcenter{
\xymatrixrowsep{1pc}
\xy\xymatrix{
4\edge[d]_\alpha\\
5}\endxy}} 
\right)^2 \oplus
\left(
\ignore{\vcenter{
\xymatrixrowsep{1pc}
\xy\xymatrix{
4\edge[d]_\beta\\
5}\endxy}} 
\right)^2.$$
Thus $\rep((\la/J)_\la) =4$, which by Observation 2.5 yields
$\lFindim\la
\le 4$. (Since the left $\la$-module $\la/\la(\alpha +\beta)$ has projective
dimension 4, we actually obtain equality.)

Moreover, we observe that $T$ has finite syzygy type, the syzygy category
of $T$ being equal to $\add(A_1,\dots,A_7)$ with $A_i\cong S_i$ for
$i=1,\dots,5$ while
\ignore{ $A_6= 
\vcenter{\xymatrixrowsep{1pc}
\xy\xymatrix{
4\edge[d]_\alpha\\
5}\endxy}$
and $A_7=
\vcenter{\xymatrixrowsep{1pc}
\xy\xymatrix{
4\edge[d]_\beta\\
5}\endxy}$.
}
Using the notation introduced above, we thus obtain $s=7$ and the integer
$7\times7$ matrix
$$\bfb= \pmatrix 0&1&0&0&0&0&0\\
0&0&1&0&0&0&0\\
0&0&0&1&0&0&0\\
0&0&0&0&2&1&1\\
0&0&0&0&0&0&0\\
0&0&0&0&1&1&1\\
0&0&0&0&1&1&1\endpmatrix.$$
The projective covers of the $A_i$ being $P_i$ for $i\le5$ and $P_6\cong
P_7\cong e_4\la$, we finally see that $p_{ii}=1$ for $i\le 5$ and $p_{64}=
p_{74}=1$, while all the other exponents $p_{il}$ are zero.

We will use Proposition 3.1 to determine the projective dimension of the
left $\la$-module 
$$M= \bigl( \la e_4\oplus \la e_5 \bigr) \bigm/ \bigl( \la(\gamma,\alpha)
+\la(\beta,\delta) +\la(\epsilon,0) \bigr).$$ 
Clearly, $\dim_K(A_i\otimes_\la M)=0$ for
$i=1,2,3$ and 
$$\dim_K(A_4\otimes_\la M)= \dim_K(e_4M/e_4JM) =1= \dim_K(A_5\otimes_\la
M).$$
Moreover, we compute that
$$\dim_K(A_6\otimes_\la M)= \dim_K\bigl( e_4M/ (\gamma\la +\delta\la
+\beta\la)M \bigr) =1,$$
and analogously, $\dim_K(A_7\otimes_\la M)=1$. In view of 3.1(c), this
yields
$$\tau(M)= \tau= (0,0,1,2,0,1,1)^{\tr}.$$
Since the last three rows of each power $\bfb^m$ are of the form
$$\pmatrix 0&0&0&0&0&0&0\\
0&0&0&0&2^m&2^m&2^m\\
0&0&0&0&2^m&2^m&2^m \endpmatrix$$
for $m\in\NN$, we deduce that $\bfb^m\tau \ne 0$ for all $m$. Thus 3.1(b)
tells us that $\pdim M=\infty$. \qed\enddefinition

\head 4. A new class of algebras of finite global repetition
index\endhead

As we will see, any classical order $\O$ of finite global dimension $d$
over a discrete valuation ring $D$ has a finite dimensional satellite
algebra $\la$ which has (left and right) global repetition index $d-1$.
This algebra, first studied by Tarsy \cite{\Tartwo} and Jategaonkar
\cite{\Jat}, and subsequently by Kirkman-Kuzmanovich \cite{\KiKu}, faithfully
reflects the homological behavior of $\O$ while being structurally less
involved.

We start with a synopsis of problems and results pertaining to classical
orders. Let $D$ be a discrete valuation ring with uniformizing parameter
$\pi$ (that is, $\pi D$ is the maximal ideal of $D$), and quotient field
$Q$. Moreover, let $\O$ be a classical order in some matrix ring $M_n(Q)$,
that is, $\O\subseteq M_n(Q)$ is a $D$-subalgebra which is finitely
generated as a $D$-module on one hand and which generates $M_n(Q)$ over $Q$
on the other. We will identify $D$ with the subring
$D\cdot1$ of $\O$.

In 1970, Tarsy \cite{\Tar} conjectured that the global dimension of any
classical order $\O\subseteq M_n(Q)$, whenever finite, is bounded above by
$n-1$. At least in the special case where $\O$ is {\it tiled}, meaning that
$\O$ contains a complete set of $n$ primitive orthogonal idempotents of
$M_n(Q)$, this conjecture seemed to have some merit. Tarsy showed in
\cite{\Tartwo} that there are only finitely many possible finite global
dimensions for tiled classical orders in $M_n(Q)$, thus
guaranteeing a finite upper bound on these dimensions. Moreover, in
\cite{\Jat} Jategaonkar proved that, up to conjugation, $M_n(Q)$ contains
only finitely many tiled classical orders of finite global dimension.
Subsequently, Tarsy's conjecture was confirmed in the triangular case (i.e.,
when $\O$ is conjugate to an order of the form
$$\pmatrix D &D &\cdots &D\\
(\pi^{\lambda_{2,1}}) &D &\cdots &D\\
\vdots &\ddots &\ddots &\vdots\\
(\pi^{\lambda_{n,1}}) &\cdots &(\pi^{\lambda_{n,n-1}}) &D\endpmatrix$$
with $\lambda_{i,j}\ge 0$) in \cite{\Jat}, and for tiled orders in $M_n(D)$
containing the ideal $M_n(\pi D)$ in \cite{\KiKu}. However,
shortly afterwards, Fujita \cite{\Fuj} constructed a class of tiled examples
$\O_n
\subseteq M_n\bigl( K((\pi)) \bigr)$ for $n\ge 6$, where $K$ is a field, such
that
$\gldim
\O_n =n$, thus refuting the conjecture even in the tiled case. Moreover, in
view of an example of Rump (\cite{\Rum}, \S7, Example 3) and a class of
examples exhibited by Jansen and Odenthal in
\cite{\JaOd}, raising `$n-1$' to `$n$' (or even `$n+1$') in Tarsy's
conjecture will not suffice to provide a valid bound either. 

The satellite algebra of $\O$ which we will consider here is the factor
$\la= \opio$. Clearly, $\la$ is a finite dimensional algebra over the residue
field $K= D/\pi D$. For the case where $\O$ is tiled and basic, Kirkman and
Kuzmanovich assembled a list of distinguished properties of $\la$,
including the following \cite{\KKprivate}:

$\bullet$ $\la$ is a split basic algebra. More precisely, $\la\cong
K\Gamma/I$ is a binomial relation algebra, where the ideal $I$ has a
generating set consisting of paths and differences of paths. Furthermore,
the quiver
$\Gamma$ and a generating set for $I$ of the described form can be
explicitly (and easily) computed from the valued quiver which was associated
to any tiled classical order by Wiedemann and Roggenkamp \cite{\WiRo}. (A
sketch of this procedure can be found ahead of Example 4.6 below.)

$\bullet$ Each of the simple $\la$-modules occurs as a composition factor
of multiplicity 1 in each of the indecomposable projective $\la$-modules.
Consequently, the multiplicities of the simples are equal in all
$\la$-modules of finite projective dimension, which, in particular, shows
that all simple modules have infinite projective dimension for $n\ge2$.

Since $\O$ is noetherian and $\pi$ is a central non-zero-divisor in $\O$
which is not a unit, the standard change-of-rings
arguments for homological dimensions (see, e.g., \cite{\Kap}) apply to the
pair $\O$ and $\la=\opio$. They yield a helpful relationship between the
little finitistic dimensions of $\O$ and $\la$, as was observed by Tarsy
and by Green, Kirkman, and Kuzmanovich:

\proclaim{Proposition 4.1} $\lfindim\O= \lfindim\la+1$ and $\rfindim\O=
\rfindim\la+1$.\endproclaim

\demo{Proof} \cite{\Tartwo}, proof of Corollary 1; \cite{\GKK}, Lemma 2.7.
\qed\enddemo

For the big finitistic dimensions, however, these arguments only
yield inequalities
$$\lFindim\O \ge\lFindim\la+1 \qquad \text{and} \qquad \rFindim\O
\ge\rFindim\la+1.$$
In case the order $\O$ has finite global dimension, or more generally, if it
has finite injective dimension on both sides, we can combine Proposition 4.1
with Proposition 1.3 to obtain coincidence of the big and little finitistic
dimensions of
$\la$, as follows.

\proclaim{Proposition 4.2} If $\idim {}_\O\O$ and $\idim \O_\O$ are both
finite, then
$$\align \lFindim\la &= \rFindim\la= \lfindim\la= \rfindim\la\\
 &= \idim {}_\la\la= \idim \la_\la = \pdim E(_\la(\la/J))=
\pdim E((\la/J)_\la)\\
 &= \lfindim\O-1 =\rfindim\O-1. \endalign$$
In particular, if $\gldim\O= d<\infty$, then $d\ge1$ and all the above
dimensions equal $d-1$.\endproclaim

\demo{Proof} Since the order $\O$ is not divisible by $\pi$ on either side,
it is not injective as a right or left $\O$-module. Therefore, $\pi$ being
a central non-zero-divisor in $\O$, it follows from \cite{\Kap}, Theorem 205
that
$$\idim {}_\la\la \le \idim {}_\O\O-1 <\infty,$$
and likewise that $\idim \la_\la <\infty$. Thus, by Proposition 1.3,
$$\align \lFindim\la &= \rFindim\la= \lfindim\la= \rfindim\la\\
 &= \idim {}_\la\la= \idim \la_\la = \pdim E(_\la(\la/J))=
\pdim E((\la/J)_\la). \endalign$$
If $d'$ is the common value of these dimensions, Proposition 4.1 shows that
$\lfindim\O= \rfindim\O= d'+1$. 

The final statement of the proposition is clear. \qed\enddemo

We shall need the following well known variant of Schanuel's Lemma, which
holds for modules over an arbitrary artinian ring: If
$0\rightarrow L\rightarrow P\rightarrow N\rightarrow 0$ is a short exact
sequence of
$\la$-modules with
$P$ projective, then $L\cong \Omega^1(N)\oplus Q$ for some projective
$\la$-module $Q$. This conclusion can be obtained from Schanuel's Lemma
together with the Krull-Schmidt Theorem, or from \cite{\AnFu}, Lemma 17.17.

\proclaim{Theorem 4.3} If $M$ is a nonzero left $\la$-module with $\pdim
{}_\O M =m<\infty$, then $m>0$ and $\Omega^{m-1}_\la(M) \cong
\Omega^{m+1}_\la(M)\oplus Q$ for some projective module $Q$. In particular,
$\rep(_\la M)\le m-1$, and if $M$ is finitely generated, then $_\la M$ has
finite syzygy type.\endproclaim

\demo{Proof} Since $M$ is torsion as an $\O$-module, $m\ge1$.

Suppose first that $m=1$, and write $M\cong P/K$ for some projective
$\O$-modules $K\subseteq P$. Since $\pi M=0$, there is an exact sequence
$$0\rightarrow K/\pi P\rightarrow P/\pi P\rightarrow M\rightarrow 0$$
of $\la$-modules with $P/\pi P$ projective, and hence $K/\pi
P\cong \Omega^1_\la(M)\oplus Q_1$ for some projective $\la$-module $Q_1$.
There is also an exact sequence 
$$0\rightarrow \pi P/\pi K\rightarrow K/\pi
K\rightarrow K/\pi P\rightarrow 0$$
of $\la$-modules with $K/\pi K$
projective, and so $\pi P/\pi K\cong \Omega^1_\la(K/\pi P)\oplus Q$ for some
projective $\la$-module $Q$. Consequently, $\pi P/\pi K\cong \Omega^2_\la(M)
\oplus Q$. Since $\pi$ is a non-zero-divisor on $P$, we have $\pi P/\pi
K\cong P/K\cong M$ as well. 

Now assume that $m>1$, and consider the exact sequence $0\rightarrow
L\rightarrow Q_0\rightarrow M\rightarrow 0$ of $\la$-modules where $Q_0$ is
a $\la$-projective cover of $M$ and $L= \Omega^1_\la(M)$. Observe that $\pdim
{}_\O\la=1$ and so $\pdim {}_\O Q_0\le 1$. Since $m\ge2$, it follows from
the long exact sequence for $\Ext_\O$ that
$\pdim {}_\O L= m-1$. By induction, $\Omega_\la^{m-2}(L)\cong
\Omega_\la^m(L)\oplus Q$ for some projective module $Q$, and therefore
 $\Omega^{m-1}_\la(M) \cong
\Omega^{m+1}_\la(M)\oplus Q$.

The final conclusions of the theorem now follow easily. \qed\enddemo

\proclaim{Corollary 4.4} If $\gldim\O =d<\infty$, then $\lglrep\la
=\rglrep\la = d-1$ and all finitely generated $\la$-modules have finite
syzygy type.
\endproclaim

\demo{Proof} Suppose that $\gldim\O =d<\infty$. Theorem 4.3 then shows that
$\lglrep\la$ and $\rglrep\la$ are bounded above by $d-1$, and that all
finitely generated $\la$-modules have finite syzygy type. To obtain the
missing inequalities, we use Proposition 4.1 and Observation 2.5 to obtain
$$d-1=
\lfindim\la
\le
\rep((\la/J)_\la)
\le \rglrep\la,$$
and similarly $d-1\le \lglrep\la$. \qed\enddemo

The next consequence of Theorem 4.3 addresses the Poincar\'e-Betti series
of a pair $(M,N)$ of $\la$-modules, i.e. the formal power series
$$\sum_{i=0}^\infty (-1)^i \dim_K \Ext^i_\la(M,N)\cdot t^i.$$

\proclaim{Corollary 4.5} If $\O$ has finite global dimension, then the
Poincar\'e-Betti series of any pair of $\la$-modules is a rational
function. \endproclaim

\demo{Proof} Combine Corollary 4.4 with Wilson's main result in
\cite{\Wil}. \qed\enddemo 

We now briefly sketch the road that leads from a basic tiled classical order
$\O$ to a presentation of $\la=\opio$ in terms of quiver and relations, as
it was communicated to us by Kirkman and Kuzmanovich \cite{\KKprivate}. We
would like to warn the reader, however, that the quivers we use are duals
of those considered in \cite{\KiKu} and \cite{\WiRo}.

After a conjugation, we may
assume that
$\O\subseteq M_n(D)$ and that the diagonal matrix units $e_{ii}$, $1\le
i\le n$, all lie in $\O$. Now consider the ordinary quiver
$\Gamma$ of the (semiperfect) ring $\O$: It consists of $n$ vertices
$1,\dots,n$ corresponding to the $n$ non-isomorphic indecomposable
projective left $\O$-modules $P_i= \O e_{ii}$, and there is an arrow
$i\rightarrow j$ provided that $P_j$ is isomorphic to a direct summand of
the projective cover of $\rad P_i$. (This is in accordance with the usual
convention, since the multiplicity of $P_j$ as a direct summand of the
projective cover of $\rad P_i$ is at most 1 in the present situation.)
We attach to any arrow $i@>\alpha>> j$ the value
$v(\alpha)=k$ where $k$ is the least nonnegative integer such that
$\pi^kP_j\subseteq \rad P_i$ (so that $D\pi^ke_{ji}$ corresponds, via
right multiplication, to $\Hom_{\O}(P_j, \rad P_i)$). Thus $D\pi^ke_{ji}=
e_{jj}\O e_{ii}$ if $i\ne j$, whereas $k=1$ if $i=j$. This
process yields the `valued quiver' $\tlgam$ of $\O$. (The order $\O$ can
actually be retrieved from $\tlgam$ -- see \cite{\WiRo}.) The value
function $v$ is extended additively from arrows to paths. Thus the value of
a path $q= \alpha_m\alpha_{m-1} \cdots\alpha_1$ in $\tlgam$, where the
$\alpha_i$ are arrows, is $v(q)= \sum_{i=1}^m v(\alpha_i)$.

We now specialize to the situation where $\tlgam$ has no loops. In that
case the quiver of $\la=\opio$ coincides with the quiver $\Gamma$ as above,
and the following set $G$ of paths and differences of paths in
$K\Gamma$ generates an ideal $I\subseteq K\Gamma$ such that $\la\cong
K\Gamma/I$. A path $q$, with starting point $i$ and endpoint $j$ say,
belongs to $G$ if and only if there exists a path $q'$ from $i$ to $j$ in
$\tlgam$ such that $v(q') <v(q)$. Moreover, given two paths $p$ and $q$
with coinciding starting points and coinciding endpoints, the difference
$p-q$ belongs to $G$ if and only if $v(p)= v(q)$.
\medskip

Since, given a tiled classical order $\O$, the algebra $\la=\opio$ is
substantially easier to handle, this latter algebra serves as an excellent
tool for the investigation of $\O$. In our first illustration of this
principle, we exhibit a classical order with differing left and right
finitistic dimensions.

\definition{Example 4.6} Let $D$ be any discrete valuation ring with
uniformizing parameter $\pi$ and quotient field $Q$, and consider the
following tiled classical order $\O\subseteq M_6(Q)$:
$$\O = \pmatrix D&D&D&D&D&D\\
(\pi) & D &(\pi) &(\pi) &D &D\\
(\pi) &(\pi) &D &D &D &D\\
(\pi^2) &(\pi) &(\pi^2) &D &(\pi) &D\\
(\pi^2) &(\pi) &(\pi) &(\pi) &D &D\\
(\pi^2) &(\pi^2) &(\pi^2) &(\pi) &(\pi) &D \endpmatrix$$
The valued quiver of $\O$ is

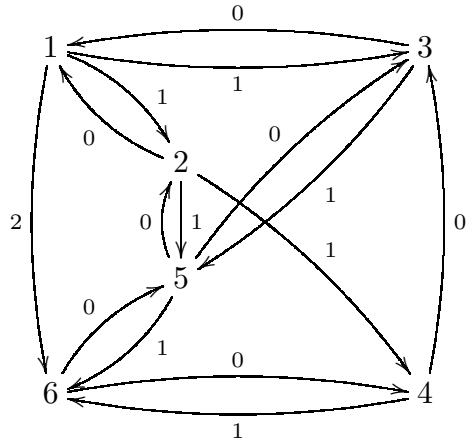
\ignore{
$$\xymatrixcolsep{2.5pc}
\xy\xymatrix{
1 \ar@/_/[ddd]_2 \ar@/^/[dr]^(0.65)1 \ar@/_/[rrr]_1 &&&3 \ar@/^/[ddll]^1
\ar@/_/[lll]_0\\
 &2 \ar@/^/[ul]^0 \ar[d]^1 \ar@/^/[ddrr]^1\\
 &5 \ar@/^/[u]^0 \ar@/^/[uurr]^0 \ar@/^/[dl]^(0.35)1\\
6 \ar@/^/[rrr]^0 \ar@/^/[ur]^0 &&&4 \ar@/^/[lll]^1 \ar@/_/[uuu]_0
}\endxy$$
}

\noindent and $\la=\opio$ has the same quiver (non-valued) and relations
$p$ and $p-q$, for suitable paths $p$ and $q$ such that the indecomposable
projective \underbar{right} $\la$-modules are as in Example 2.8. Combining
our computations in the latter example with Proposition 4.1, we thus obtain
$\lfindim \O= 4$ and $\rfindim \O= 1$. \qed\enddefinition

In our final example, we construct a tiled classical order $\O$ of infinite
global dimension such that $\la= \opio$ is Gorenstein, i.e., such that
$\idim {}_\la\la =\idim \la_\la <\infty$. The fact that $\gldim \O= \infty$
will be established by pinpointing a left $\la$-module which does not
satisfy the conclusion of Theorem 4.3. 

\definition{Example 4.7} Again, let $D$ be any discrete valuation ring with
uniformizing parameter $\pi$ and quotient field $Q$, and let
$\O\subseteq M_6(Q)$ be the tiled classical order which has valued quiver

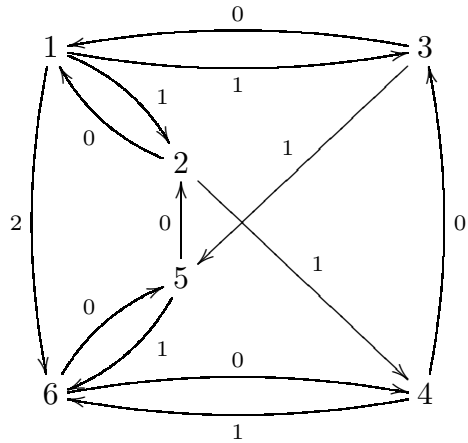
\ignore{
$$\xymatrixcolsep{2.5pc}
\xy\xymatrix{
1 \ar@/_/[ddd]_2 \ar@/^/[dr]^(0.65)1 \ar@/_/[rrr]_1 &&&3 \ar[ddll]_1
\ar@/_/[lll]_0\\
 &2 \ar@/^/[ul]^0 \ar[ddrr]^1\\
 &5 \ar[u]^0 \ar@/^/[dl]^(0.35)1\\
6 \ar@/^/[rrr]^0 \ar@/^/[ur]^0 &&&4 \ar@/^/[lll]^1 \ar@/_/[uuu]_0
}\endxy$$
}

\noindent Then the quiver $\Gamma$ of $\la=\opio$ is the non-valued version
of the one shown and $\la$ is of the form $K\Gamma/I$, where $I$ is
generated by paths and differences of paths, such that the indecomposable
projective \underbar{left} $\la$-modules have the following graphs:

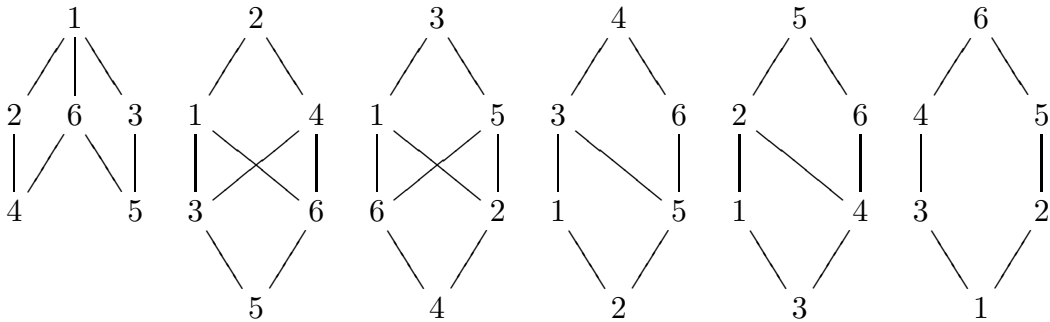
\ignore{
$$\xymatrixcolsep{0.67pc}\xymatrixrowsep{1.5pc}
\xy\xymatrix{
 &1 \edge[dl] \edge[d] \edge[dr] &&&2 \edge[dl] \edge[dr] &&&3 \edge[dl]
\edge[dr] &&&4 \edge[dl] \edge[dr] &&&5 \edge[dl] \edge[dr] &&&6
\edge[dl] \edge[dr]\\
2 \edge[d] &6 \edge[dl] \edge[dr] &3 \edge[d] &1 \edge[d] \edge[drr] &&4
\edge[dll] \edge[d] &1 \edge[d] \edge[drr] &&5 \edge[dll]
\edge[d] &3 \edge[d]
\edge[drr] &&6 \edge[d] &2 \edge[d] \edge[drr] &&6 \edge[d] &4 \edge[d]
&&5 \edge[d]\\
4 &&5 &3 \edge[dr] &&6 \edge[dl] &6 \edge[dr] &&2 \edge[dl] &1 \edge[dr]
&&5 \edge[dl] &1 \edge[dr] &&4 \edge[dl] &3 \edge[dr] &&2 \edge[dl]\\
 &&&&5 &&&4 &&&2 &&&3 &&&1
}\endxy$$
}

\noindent In particular, the only indecomposable injective object in
$\lamod$ which fails to be projective is the injective envelope $E_6$ of
$\la e_6/Je_6$, namely the unique left $\la$-module with graph

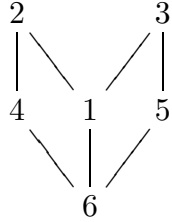
\ignore{
$$\xymatrixcolsep{1pc}\xymatrixrowsep{1.5pc}
\xy\xymatrix{
2 \edge[d] \edge[dr] &&3 \edge[dl] \edge[d]\\
4 \edge[dr] &1 \edge[d] &5 \edge[dl]\\
 &6
}\endxy$$
}

\noindent One checks that $\Omega^1(E_6)\cong \la e_1$, whence $\pdim E_6
=1$. Thus the minimal injective cogenerator for $\lamod$ has projective
dimension 1. Analogously, one computes the projective dimension of the
minimal injective right cogenerator to be 1, which yields
$\idim {}_\la\la= \idim \la_\la =1$ by duality. Furthermore, $\lfindim\la
=\rfindim\la =1$ by Proposition 1.3. Consequently,
Proposition 4.1 tells us that $\lfindim\O =\rfindim\O =2$.

In order to check that $\gldim\O =\infty$, it thus suffices to procure a
left $\la$-module $M$ which violates the conclusion of Theorem 4.3 for
$m=2$. One can readily verify that each of the simple left $\la$-modules
qualifies for that purpose. We pick one with a `slim' projective
resolution, namely $M= \la e_6/Je_6$. Then the first three syzygies
$\Omega^1(M)$, $\Omega^2(M)$, $\Omega^3(M)$ have the following graphs,
respectively:

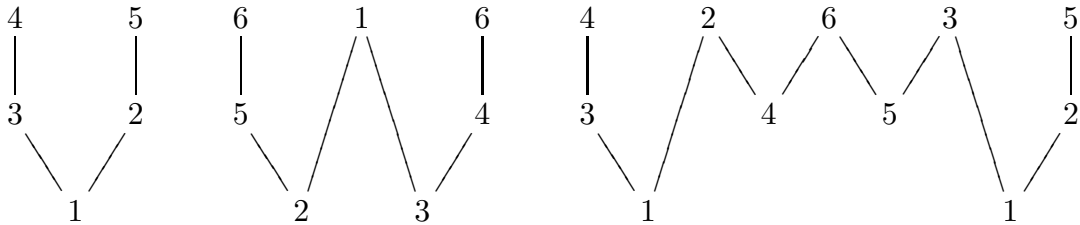
\ignore{
$$\xymatrixcolsep{0.67pc}\xymatrixrowsep{1.5pc}
\xy\xymatrix{
4 \edge[d] &&5 \edge[d] &&6 \edge[d] &&1 \edge[ddl] \edge[ddr] &&6
\edge[d] &&4 \edge[d] &&2 \edge[ddl] \edge[dr] &&6 \edge[dl] \edge[dr]
&&3 \edge[dl] \edge[ddr] &&5 \edge[d]\\
3 \edge[dr] &&2 \edge[dl] &&5 \edge[dr] &&&&4 \edge[dl] &&3 \edge[dr]
&&&4 &&5 &&&2 \edge[dl]\\
 &1 &&&&2 &&3 &&&&1 &&&&&&1
 }\endxy$$
}

\noindent Clearly, these syzygies are indecomposable and $\Omega^1(M)
\not\cong \Omega^3(M)$, which shows that indeed $\pdim_\O M >2$. (The
repetition index of $M$ as a $\la$-module is, in fact, infinite.)
\qed\enddefinition

\enddocument